\documentstyle{article}
\newtheorem{rem}{Remark}

\begin{document}
\begin{center}
\Large{\bf  RIEMANN EXTENSIONS IN THEORY OF THE FIRST\\[2mm] ORDER
SYSTEMS OF DIFFERENTIAL EQUATIONS}\vspace{4mm}\normalsize
\end{center}
 \begin{center}
\Large{\bf Valery Dryuma}\vspace{4mm}\normalsize
\end{center}
\begin{center}
{\bf Institute of Mathematics and Informatics AS Moldova, Kishinev}\vspace{4mm}\normalsize
\end{center}
\begin{center}
{\bf  E-mail: valery@dryuma.com;\quad cainar@mail.md}\vspace{4mm}\normalsize
\end{center}
\begin{center}
{\bf  Abstract}\vspace{4mm}\normalsize
\end{center}

   The properties of the Riemann extensions of nonriemannian spaces defined by the first order
    systems of differential equations are considered.

\section{Introduction}

   The first order polynomial  systems of differential equations
\begin{equation} \label{dryuma:eq1}
\frac{d x^{i}}{ds}=a^{i}_{j}x^{j}+b^{i}_{ j k}x^{j}x^{k}
\end{equation}
depending from the parameters $a,b$  play an important role in
various branches of modern mathematics and its applications.

   However even in the case of the  system of two equations
\begin{equation} \label{dryuma:eq2}
\frac{d x}{ds}=k+ax+by+cx^2+exy+fy^2,\quad \frac{d
y}{ds}=l+mx+ny+px^2+qxy+ry^2
\end{equation}
there are many unsolved problems.

    The most famous of them is the 16-th Hilbert problem about the
quantity and position of the limit cycles in such type of the
systems.

  The spatial first order system of differential equations
  \begin{equation} \label{dryuma:eq3}
\frac{d x}{ds}=P(x,y,z),\quad \frac{d y}{ds}=Q(x,y,z),\quad
\frac{d z}{ds}=R(x,y,z)
\end{equation}
with the functions  $P,Q,R$ polynomial on variables $x,y,z$
 are still more complicated object for the
studying of their properties.

As example the studying of simplest spatial systems of equations such as the Lorenz
equations
\begin{equation} \label{dryuma:eq4}
\frac{d x}{ds}=\sigma(y-x),\quad \frac{d y}{ds}=rx-y-xz,\quad
\frac{d z}{ds}=xy-bz
\end{equation}
or the R\"ossler system
\begin{equation} \label{dryuma:eq5}
\frac{d x}{ds}=-y-z,\quad \frac{d y}{ds}=x+ay,\quad \frac{d
z}{ds}=bx-cz+xz
\end{equation}
having chaotic behavior at some values of parameters represent the
difficult task.

   The systems of the first order differential equations are not suitable
   object of consideration from the usually point of Riemann geometry.

   The systems of the second order differential equations in form
\begin{equation} \label{dryuma:eq6}
\frac{d^2 x^i}{ds^2}+\Pi^i_{k j}(x)\frac{d x^k}{ds}\frac{d
x^j}{ds}=0
\end{equation}
are best suited to do that.

   They can be considered as geodesics of the affinely connected space
   $M$ in local coordinates $x^k$. The values $\Pi^i_{jk}=\Pi^i_{kj}$ are the coefficients of
   affine connections on $M$.

    With the help of such coefficients can
   be constructed curvature tensor and others geometrical objects
   defined on variety $M$.

    There are many possibilities to present a given system of the first order of equations
    in the form of (\ref{dryuma:eq6}).

    One of them is a following presentation.

    For the system
 \begin{equation} \label{dryuma:eq7}
\frac{d x}{ds}=P(x,y),\quad \frac{d y}{ds}=Q(x,y)
\end{equation}
after differentiation with respect to parameter $s$  we get the
second order system of differential equations of the form
(\ref{dryuma:eq6})
\begin{equation}\label{dryuma:eq8}
\frac{d^2 x}{ds^2}=\frac{1}{P}(P_x \frac{d x}{ds}+P_y \frac{d
y}{ds})\frac{d x}{ds}$$  ,
 $$
 \frac{d^2 y}{ds^2}=\frac{1}{Q}(Q_x
\frac{d x}{ds}+Q_y \frac{d y}{ds})\frac{d y}{ds}.
\end{equation}

    Such type of the system contains the integral curves of
    the system (\ref{dryuma:eq7}) as part of its solutions and can be considered as
    the equations of geodesics of two dimensional space $M^2(x,y)$ equipped by affine
    connections with coefficients
\begin{equation} \label{dryuma:eq9}
\Pi^1_{11}=-\frac{P_x}{P},\quad \Pi^1_{12}=-\frac{P_y}{2P},\quad
\Pi^2_{12}=-\frac{Q_x}{2Q},\quad \Pi^2_{22}=-\frac{Q_y}{Q}.
\end{equation}

    It is apparent that the properties of the system
    (\ref{dryuma:eq7}) have an influence on geometry of the variety $M^2(x,y)$.

    Remark that the system (\ref{dryuma:eq8}) is equivalent the
    second order differential equation
\[
\frac{d^2 y}{dx^2}=\ln(Q/P)_y\left(\frac{d
y}{dx}\right)^2+\ln(Q/P)_x\frac{d y}{dx}
\]
which has the solution in form
\[
\frac{d y}{dx}=\frac{Q}{P}.
\]

    In this connection we can use the second order system of differential
    equation
\begin{equation}\label{dryuma:eq10}
 \frac{d^2 x}{ds^2}=-\ln(Q/P)_x\left(\frac{d
x}{ds}\right)^2,\quad \frac{d^2 y}{ds^2}=\ln(Q/P)_y\left(\frac{d
y}{ds}\right)^2
\end{equation}
 for the studying of the properties of the first order system of equations (\ref{dryuma:eq7}).

    By analogy can be written the spatial system of the first
    order differential equations.

\section{The Riemann extension of affinely connected space}

      For the studying of the geometry of the equations like  (\ref{dryuma:eq8})
      we apply the notion of the Riemann extension of nonriemannian  space which was used
      earlier in  \cite{dryuma1:dryuma,dryuma2:dryuma,dryuma3:dryuma}.

     Remind basic properties of this construction.

     With help of the coefficients of affine connection of a given n-dimensional space
      can be introduced  2n-dimensional
     Riemann space $D^{2n}$ in local coordinates $(x^i,\Psi_i)$ having the metric of form

\begin{equation} \label{dryuma:eq11}
{^{2n}}ds^2=-2\Pi^k_{ij}(x^l)\Psi_k dx^i dx^j+2d \Psi_k dx^k
\end{equation}

\noindent where $\Psi_{k}$ are the additional coordinates.

The important property of such type metric is that the geodesic
 equations of metric (\ref{dryuma:eq11})  decomposes into two parts
\begin{equation} \label{dryuma:eq12}
\ddot x^k +\Pi^k_{ij}\dot x^i \dot x^j=0,
\end{equation}
and
\begin{equation} \label{dryuma:eq13}
\frac{\delta^2 \Psi_k}{ds^2}+R^l_{kji}\dot x^j \dot x^i \Psi_l=0,
\end{equation}
where
\[
\frac{\delta \Psi_k}{ds}=\frac{d
\Psi_k}{ds}-\Pi^l_{jk}\Psi_l\frac{d x^j}{ds}
\]
and $R^l_{kji}$ are the curvature tensor of n-dimensional space
with a given affine connection.

 The first part (\ref{dryuma:eq12}) of the full system
is the system of equations for geodesic of basic space with local
coordinates $x^i$ and it do not contains the supplementary
coordinates $\Psi_k$.

 The second part (\ref{dryuma:eq13}) of the system  has the form
of linear $N\times N$ matrix system of second order ODE's for
supplementary  coordinates $\Psi_k$
\begin{equation} \label{dryuma:eq14}
\frac{d^2 \vec \Psi}{ds^2}+A(s)\frac{d \vec \Psi}{ds}+B(s)\vec
\Psi=0.
\end{equation}

   Remark that the full system of geodesics has the first integral
\begin{equation} \label{dryuma:eq15}
-2\Pi^k_{ij}(x^l)\Psi_k \frac{dx^i}{ds}\frac{dx^j}{ds}+2\frac{d
\Psi_k}{ds}\frac{dx^k}{ds}=\nu
\end{equation}
which is equivalent to the relation
\begin{equation} \label{dryuma:eq16}
2\Psi_k\frac{dx^k}{ds}=\nu s+\mu
\end{equation}
where $\mu, \nu$ are parameters.

   It is important to note that the geometry of extended space
connects with geometry of basic space. For example the property of
the space to be Ricci-flat $R_{ij}=0$ or symmetrical
$R_{ijkl;m}=0$ keeps also for the extended space.

    It is important to note that for extended  space having the metric (\ref{dryuma:eq11})
    all scalar curvature invariants are vanished.

    As consequence the properties of linear system of
equation (\ref{dryuma:eq13}-\ref{dryuma:eq14}) depending from the
the invariants of  $N\times N$ matrix-function
\[
E=B-\frac{1}{2}\frac{d A}{ds}-\frac{1}{4}A^2
\]
under change of the coordinates $\Psi_k$ can be of used for that.

 The first applications the notion of extended spaces for the studying of nonlinear second order differential
 equations connected with nonlinear dynamical systems have been considered in the works of author
\cite{dryuma1:dryuma,dryuma2:dryuma,dryuma3:dryuma}.

\section{The geometry of planar system of equations in form (\ref{dryuma:eq8}) }

    We shall consider from geometrical point  the properties of planar
    systems of the equations (\ref{dryuma:eq7}).

    To do this we use the second order system (\ref{dryuma:eq8})
    having the solutions of the system (\ref{dryuma:eq7}) as the part
    of their own solutions.

    Now the system (\ref{dryuma:eq8}) can be considered as the
    geodesics of two dimensional space with affine connections and
    we extend this space up to the four-dimensional space by the introducing two additional
    coordinates $\psi_1=z,\psi_2=t$.

    The Riemann metric of the four-dimensional extended space is
    defined by
\begin{equation} \label{dryuma:eq17}
{^4}ds^2 = 2 z \frac{P_x}{P} dx^2
+(2z\frac{P_y}{P}+t\frac{Q_x}{Q})dxdy+2t\frac{Q_y}{Q}dy^2+2dxdz+2dydt.
\end{equation}

    Let us consider the basic geometric characteristics of a given
    metric.

    The nonzero components of the Ricci tensor are
\[
R_{11}=\]\[=1/2\,{\frac {-3\,\left ({\frac {\partial }{\partial
x}}Q(x,y)\right )^ {2}P(x,y)+2\,\left ({\frac {\partial
^{2}}{\partial {x}^{2}}}Q(x,y) \right )P(x,y)Q(x,y)+2\,\left
({\frac {\partial }{\partial x}}P(x,y) \right )Q(x,y){\frac
{\partial }{\partial x}}Q(x,y)}{P(x,y)\left (Q(x, y)\right
)^{2}}},
\]
\[
R_{12}=-1/2\,{\frac {-\left ({\frac {\partial ^{2}}{\partial
x\partial y}}P(x ,y)\right )P(x,y)\left (Q(x,y)\right )^{2}-\left
({\frac {\partial }{
\partial x}}Q(x,y)\right )\left ({\frac {\partial }{\partial y}}P(x,y)
\right )P(x,y)Q(x,y)}{\left (P(x,y)\right )^{2}\left (Q(x,y)\right
)^{ 2}}}-\]\[-1/2\,{\frac {\left ({\frac {\partial }{\partial
x}}P(x,y)\right ) \left ({\frac {\partial }{\partial
y}}P(x,y)\right )\left (Q(x,y) \right )^{2}-\left ({\frac
{\partial ^{2}}{\partial x\partial y}}Q(x,y )\right )\left
(P(x,y)\right )^{2}Q(x,y)}{\left (P(x,y)\right )^{2} \left
(Q(x,y)\right )^{2}}}-\]\[-1/2\,{\frac {\left ({\frac {\partial }{
\partial x}}Q(x,y)\right ){\frac {\partial }{\partial y}}Q(x,y)}{
\left (Q(x,y)\right )^{2}}},
\]
\[
R_{22}=\]\[=1/2\,{\frac {-3\,\left ({\frac {\partial }{\partial
y}}P(x,y)\right )^ {2}Q(x,y)+2\,\left ({\frac {\partial
^{2}}{\partial {y}^{2}}}P(x,y) \right )P(x,y)Q(x,y)+2\,\left
({\frac {\partial }{\partial y}}Q(x,y) \right )P(x,y){\frac
{\partial }{\partial y}}P(x,y)}{Q(x,y)\left (P(x, y)\right
)^{2}}}.
\]

   As was mentioned above all scalar invariants constructed from the
   curvature tensor and its covariant derivatives
\[
p=R_{ij}R^{ij}=0,\quad q=R_{ijkl}R^{ijkl}=0 ...
\]
are vanish.

     The full system of geodesics of the metric (\ref{dryuma:eq17}) consists
     from the two groups of equations
\[
{\frac {d^{2}}{d{s}^{2}}}x \left( s \right) -{\frac { \left( {\frac {
\partial }{\partial x}}P \left( x,y \right)  \right)  \left( {\frac {d
}{ds}}x \left( s \right)  \right) ^{2}}{P \left( x,y \right) }}-{ \frac { \left(
{\frac {\partial }{\partial y}}P \left( x,y \right)
 \right)  \left( {\frac {d}{ds}}x \left( s \right)  \right) {\frac {d}
{ds}}y \left( s \right) }{P \left( x,y \right) }}=0,
\]
\[
{\frac {d^{2}}{d{s}^{2}}}x \left( s \right) -{\frac { \left( {\frac {
\partial }{\partial y}}Q \left( x,y \right)  \right)  \left( {\frac {d
}{ds}}y \left( s \right)  \right) ^{2}}{Q \left( x,y \right) }}-{ \frac { \left(
{\frac {\partial }{\partial x}}Q \left( x,y \right)
 \right)  \left( {\frac {d}{ds}}x \left( s \right)  \right) {\frac {d}
{ds}}y \left( s \right) }{Q \left( x,y \right) }}=0,
\]
and
\[
{\frac {d^{2}}{d{s}^{2}}}z \left( s \right) +{\frac { \left(  \left( P
 \left( x,y \right)  \right) ^{3} \left( {\frac {\partial }{\partial x
}}Q \left( x,y \right)  \right) ^{2}+Q \left( x,y \right)  \left( { \frac {\partial
}{\partial x}}Q \left( x,y \right)  \right)  \left( { \frac {\partial }{\partial
y}}Q \left( x,y \right)  \right)  \left( P
 \left( x,y \right)  \right) ^{2} \right) t}{Q \left( x,y \right)
 \left( P \left( x,y \right)  \right) ^{2}}}+\]\[+{\frac { \left(  \left( P
 \left( x,y \right)  \right) ^{2}Q \left( x,y \right)  \left( {\frac {
\partial }{\partial x}}Q \left( x,y \right)  \right) {\frac {\partial
}{\partial y}}P \left( x,y \right) +2\, \left( P \left( x,y \right)
 \right) ^{2} \left( Q \left( x,y \right)  \right) ^{2}{\frac {
\partial ^{2}}{\partial x\partial y}}P \left( x,y \right)  \right) z}{
Q \left( x,y \right)  \left( P \left( x,y \right)  \right) ^{2}}}+\]\[+{ \frac {
\left( - \left( Q \left( x,y \right)  \right) ^{3} \left( { \frac {\partial
}{\partial y}}P \left( x,y \right)  \right) ^{2}+
 \left( P \left( x,y \right)  \right) ^{3}Q \left( x,y \right) {\frac
{\partial ^{2}}{\partial {x}^{2}}}P \left( x,y \right) + \left( P
 \left( x,y \right)  \right) ^{2}Q \left( x,y \right)  \left( {\frac {
\partial }{\partial x}}P \left( x,y \right)  \right) ^{2} \right) z}{Q
 \left( x,y \right)  \left( P \left( x,y \right)  \right) ^{2}}}+\]\[+{
\frac { \left(  \left( Q \left( x,y \right)  \right) ^{3} \left( { \frac {\partial
^{2}}{\partial {y}^{2}}}P \left( x,y \right)  \right) P \left( x,y \right) + \left(
Q \left( x,y \right)  \right) ^{2}
 \left( {\frac {\partial }{\partial y}}Q \left( x,y \right)  \right) P
 \left( x,y \right) {\frac {\partial }{\partial y}}P \left( x,y
 \right)  \right) z}{Q \left( x,y \right)  \left( P \left( x,y
 \right)  \right) ^{2}}}+\]\[+ \left( {\frac {\partial }{\partial x}}Q
 \left( x,y \right)  \right) {\frac {d}{ds}}t \left( s \right) +{
\frac { \left( 2\, \left( {\frac {\partial }{\partial x}}P \left( x,y
 \right)  \right)  \left( P \left( x,y \right)  \right) ^{2}Q \left( x
,y \right) + \left( {\frac {\partial }{\partial y}}P \left( x,y
 \right)  \right)  \left( Q \left( x,y \right)  \right) ^{2}P \left( x
,y \right)  \right) {\frac {d}{ds}}z \left( s \right) }{Q \left( x,y
 \right)  \left( P \left( x,y \right)  \right) ^{2}}}
=\]\[=0,
\]
\\[1mm]
\[
{\frac {d^{2}}{d{s}^{2}}}t \left( s \right) +{\frac { \left(  \left( P
 \left( x,y \right)  \right) ^{2} \left( {\frac {\partial }{\partial x
}}P \left( x,y \right)  \right) Q \left( x,y \right) {\frac {\partial }{\partial
x}}Q \left( x,y \right) + \left( P \left( x,y \right)
 \right) ^{3} \left( {\frac {\partial ^{2}}{\partial {x}^{2}}}Q
 \left( x,y \right)  \right) Q \left( x,y \right)  \right) t}{ \left(
Q \left( x,y \right)  \right) ^{2}P \left( x,y \right) }}+\]\[+{\frac {
 \left( - \left( P \left( x,y \right)  \right) ^{3} \left( {\frac {
\partial }{\partial x}}Q \left( x,y \right)  \right) ^{2}+ \left( Q
 \left( x,y \right)  \right) ^{3}P \left( x,y \right) {\frac {
\partial ^{2}}{\partial {y}^{2}}}Q \left( x,y \right) +2\, \left( P
 \left( x,y \right)  \right) ^{2} \left( Q \left( x,y \right)
 \right) ^{2}{\frac {\partial ^{2}}{\partial x\partial y}}Q \left( x,y
 \right)  \right) t}{ \left( Q \left( x,y \right)  \right) ^{2}P
 \left( x,y \right) }}+\]\[+{\frac { \left( P \left( x,y \right)  \left( Q
 \left( x,y \right)  \right) ^{2} \left( {\frac {\partial }{\partial y
}}P \left( x,y \right)  \right) {\frac {\partial }{\partial x}}Q
 \left( x,y \right) + \left( Q \left( x,y \right)  \right) ^{2}P
 \left( x,y \right)  \left( {\frac {\partial }{\partial y}}Q \left( x,
y \right)  \right) ^{2} \right) t}{ \left( Q \left( x,y \right)
 \right) ^{2}P \left( x,y \right) }}+\]\[+{\frac { \left(  \left( Q \left(
x,y \right)  \right) ^{3} \left( {\frac {\partial }{\partial y}}P
 \left( x,y \right)  \right) ^{2}+P \left( x,y \right)  \left( {\frac
{\partial }{\partial x}}P \left( x,y \right)  \right)  \left( {\frac {
\partial }{\partial y}}P \left( x,y \right)  \right)  \left( Q \left(
x,y \right)  \right) ^{2} \right) z}{ \left( Q \left( x,y \right)
 \right) ^{2}P \left( x,y \right) }}+\]\[+{\frac { \left( 2\, \left( Q
 \left( x,y \right)  \right) ^{2} \left( {\frac {\partial }{\partial y
}}Q \left( x,y \right)  \right) P \left( x,y \right) +Q \left( x,y
 \right)  \left( P \left( x,y \right)  \right) ^{2}{\frac {\partial }{
\partial x}}Q \left( x,y \right)  \right) {\frac {d}{ds}}t \left( s
 \right) }{ \left( Q \left( x,y \right)  \right) ^{2}P \left( x,y
 \right) }}+\]\[+ \left( {\frac {\partial }{\partial y}}P \left( x,y
 \right)  \right) {\frac {d}{ds}}z \left( s \right)
=0.
\]

  Last of the two equations looks as the system of $2\times2$ linear matrix equations
   with the coefficients $A_{ij},B_{ij}$ depending from the solutions of the
   first system of equations ( this resemble the situation with the soliton theory!).
\begin{equation} \label{dryuma:eq18}
\frac{d^2 z}{ds^2}+A_{11}\frac{d z}{ds}+A_{12}\frac{d t}{ds}+B_{11}z +B_{12}t=0,
$$
$$ \frac{d^2 t}{ds^2}+A_{21}\frac{d z}{ds}+A_{22}\frac{d t}{ds}+B_{21}z +B_{22}t=0.
\end{equation}

     Using the first integral of geodesics of the metric (\ref{dryuma:eq16})
\[
z\frac{d x}{ds}+t\frac{d y}{ds}=\mu s/2+\nu
\]
we can get from the system (\ref{dryuma:eq18}) two independent
linear second order differential equations for the each variable
$z(s)$ and $t(s)$
\[
\frac{d^2 z}{ds^2}+A(s)\frac{d z}{ds}+B(s)z=0,
\]
and
\[
\frac{d^2 t}{ds^2}+C(s)\frac{d t}{ds}+E(s)t=0.
\]

    Let us consider some examples.

    1. The system of equations
\[
\frac{d x}{ds}=-y,\quad \frac{d y}{ds}=x
\]
corresponds the second order system of equations
\[
{\frac {d^{2}}{d{s}^{2}}}x \left( s \right) -{\frac { \left( {\frac {d }{ds}}x
\left( s \right)  \right) {\frac {d}{ds}}y \left( s \right) }{ y \left( s \right)
}}=0,
\]
\[
{\frac {d^{2}}{d{s}^{2}}}y \left( s \right) -{\frac { \left( {\frac {d }{ds}}x
\left( s \right)  \right) {\frac {d}{ds}}y \left( s \right) }{ x \left( s \right)
}}=0.
\]

   From the geodesics of extended space we get the equations for coordinate $z$
\[
{\frac {d^{2}}{d{s}^{2}}}z \left( s \right) -{\frac { \left( x \left( s \right)
\right) ^{2}z \left( s \right) }{ \left( y \left( s
 \right)  \right) ^{2}}}+1/2\,{\frac { \left( 2\,y \left( s \right)
 \left( x \left( s \right)  \right) ^{2}+2\, \left( y \left( s
 \right)  \right) ^{3} \right) {\frac {d}{ds}}z \left( s \right) }{x
 \left( s \right)  \left( y \left( s \right)  \right) ^{2}}}+1/2\,{
\frac {\mu}{x \left( s \right) }}=0
\]
and
\[
{\frac {d^{2}}{d{s}^{2}}}t \left( s \right) -{\frac { \left( y \left( s \right)
\right) ^{2}t \left( s \right) }{ \left( x \left( s
 \right)  \right) ^{2}}}-1/2\,{\frac { \left( 2\, \left( x \left( s
 \right)  \right) ^{3}+2\,x \left( s \right)  \left( y \left( s
 \right)  \right) ^{2} \right) {\frac {d}{ds}}t \left( s \right) }{y
 \left( s \right)  \left( x \left( s \right)  \right) ^{2}}}+1/2\,{
\frac {\mu}{y \left( s \right) }}=0
\]
for coordinate $t$.

\section{The geometry of planar system of equations in form (\ref{dryuma:eq10})}

   The system (\ref{dryuma:eq10}) is more simple object for the studying of the properties
   of the system (\ref{dryuma:eq7}) from geometrical point of
   view.

   Really, the metric of extended space defined by the system (\ref{dryuma:eq10})
can be written in form
\begin{equation} \label{dryuma:eq19}
{^4}ds^2 =- 2 z \frac{\partial(K(x,y)}{\partial x} dx^2
+2t\frac{\partial(K(x,y)}{\partial x}dy^2+2dxdz+2dydt
\end{equation}
where $z, t$ are the supplementary coordinates and the function
$K(x,y)=\ln(Q/P)$ is determined from the relations
\[
\frac{d y}{dx}=\exp(K(x,y))
\]
or
\[
\frac{d y}{ds}=\exp(K(x,y))\frac{d x}{ds}.
\]

   Contrary to the case (\ref{dryuma:eq17}) the Riemannian space with the metric (\ref{dryuma:eq19}) is a
   Ricci-flat
\[
R_{ij}=0.
\]

    Its geodesics are defined by the system of equations
\[
\frac{d^2 x}{ds^2}+\frac{\partial K(x,y)}{\partial x}\left(\frac{d
x}{ds}\right)^2=0,
\]
\[
\frac{d^2 y}{ds^2}-\frac{\partial K(x,y)}{\partial y}\left(\frac{d
y}{ds}\right)^2=0.
\]

\section{Rigorous approach at the geometry of the planar systems}

    We consider the system of paths of two-dimensional
     space $S_2$ in form

\begin{equation} \label{dryuma:eq20}
\ddot x +\Pi^1_{11}(\dot x)^2+2\Pi^1_{12}\dot x \dot
y+\Pi^1_{22}(\dot y)^2=0, $$ $$ \ddot y +\Pi^2_{11}(\dot
x)^2+2\Pi^2_{12}\dot x \dot y+\Pi^2_{22}(\dot y)^2=0,
\end{equation}
where coefficients $\Pi^k_{ij}=\Pi^k_{ji}$.

    The Riemann extension of the space $S_2$ is determined by the metric
\begin{equation} \label{dryuma:eq21}
{^4}ds^2=-2 z \Pi^1_{11}d x^2-2 t \Pi^2_{11}d x^2-4 z \Pi^1_{12}d
x d y-4 t \Pi^2_{12}dx d y-2 z \Pi^1_{22}d y^2-2 t \Pi^2_{22}d
y^2.
\end{equation}

    A necessary condition for the equations (\ref{dryuma:eq20}) to admit a  first
    integral
\[
a_i(x,y)\frac{d x^i}{d s}=const
\]
is
\[
a_{i;j}+a_{j;i}=0,
\]
where
\[
a_{i;j}=\frac{\partial a_i}{\partial x^j}-a_k \Gamma^k_{ij},
\]
and  $\Gamma^k_{ij}$ are the Christoffel symbols of the metric
(\ref{dryuma:eq21}) .

    We apply this conditions for determination  of the coefficients of equations
    $\Gamma^k_{ij}$ using the  vector $a_i$ in
    form
\[
a_i=[Q(x,y),-P(x,y),0,0].
\]

 This means that the first order of equation
\[
\frac{d y}{d x}=\frac{Q(x,y)}{P(x,y)}
\]
or
\[
Q(x,y)d x-P(x,y)d y=0
\]
is an integral of the paths equations.

     Using such condition it is possible to find only three
     coefficients of affine connections $\Gamma^k_{ij}$.

     For determination of the rest coefficients
we use yet another the first order equation
\begin{equation} \label{dryuma:eq22}
\frac{d y}{d x}=-\frac{y(y-1)}{x(x-1)} ,
\end{equation}
with the first integral
\[
y(x)=\frac{C(x-1)}{x-C}.
\]

   The equation (\ref{dryuma:eq22}) plays  an important role in theory of
   of the planar first order system of equations (\ref{dryuma:eq2}) (\cite{dryuma4:dryuma}).

   In result the coefficients $\Pi^k_{ij}$ of the paths equation are defined by
\[
\Pi^1_{11}={\frac {\left ({\frac {\partial }{\partial x}}Q(x,y)\right )x\left (x-
1\right )}{{y}^{2}P(x,y)-yP(x,y)+Q(x,y){x}^{2}-Q(x,y)x}},
\]
\[
\Pi^1_{12}=\]\[=1/2\,{\frac {2\,yP(x,y)-2\,P(x,y)+2\,xP(x,y)+\left ({\frac {\partial
} {\partial y}}Q(x,y)\right ){x}^{2}-\left ({\frac {\partial }{\partial
y}}Q(x,y)\right )x-\left ({\frac {\partial }{\partial x}}P(x,y)\right ){x}^{2}+\left
({\frac {\partial }{\partial x}}P(x,y)\right )x}{{y}^{2
}P(x,y)-yP(x,y)+Q(x,y){x}^{2}-Q(x,y)x}},
\]
\[
\Pi^1_{22}=-{\frac {\left ({\frac {\partial }{\partial
y}}P(x,y)\right )x\left (x -1\right
)}{{y}^{2}P(x,y)-yP(x,y)+Q(x,y){x}^{2}-Q(x,y)x}}
\]
and corresponding expressions for the coefficients $\Pi^2_{11},
\Pi^2_{12},\Pi^2_{22}$.

   In result the geodesics take the form
\begin{equation}\label{dryuma:eq23}
 {\frac {d^{2}}{d{s}^{2}}}x(s)+{\frac {\left ({\frac
{\partial }{
\partial x}}Q(x,y)\right )x\left (x-1\right )\left ({\frac {d}{ds}}x(s
)\right
)^{2}}{{y}^{2}P(x,y)-yP(x,y)+Q(x,y){x}^{2}-Q(x,y)x}}+\]\[+{\frac {
\left (2\,yP(x,y)-2\,P(x,y)+2\,xP(x,y)+\left ({\frac {\partial }{
\partial y}}Q\right ){x}^{2}-\left ({\frac {\partial }{\partial y
}}Q\right )x-\left ({\frac {\partial }{\partial x}}P\right )
{x}^{2}+\left ({\frac {\partial }{\partial x}}P\right )x\right )
\left ({\frac {d}{ds}}x(s)\right ){\frac {d}{ds}}y(s)}{{y}^{2}P-y
P+Q{x}^{2}-Qx}}-\]\[-{\frac {\left ({\frac {\partial }{
\partial y}}P\right )x\left (x-1\right )\left ({\frac {d}{ds}}y(s
)\right )^{2}}{{y}^{2}P-yP+Q{x}^{2}-Qx}}=0,
\]
\[
{\frac {d^{2}}{d{s}^{2}}}y(s)-{\frac {\left ({\frac {\partial }{
\partial x}}Q(x,y)\right )y\left (y-1\right )\left ({\frac {d}{ds}}x(s
)\right
)^{2}}{{y}^{2}P(x,y)-yP(x,y)+Q(x,y){x}^{2}-Q(x,y)x}}+\]\[+{\frac {
\left (-\left ({\frac {\partial }{\partial y}}Q\right ){y}^{2}+
\left ({\frac {\partial }{\partial y}}Q\right
)y+2\,Qy-2\,Q+2\,Qx+{y}^{2}{\frac {\partial }{\partial
x}}P-y{\frac {
\partial }{\partial x}}P\right )\left ({\frac {d}{ds}}x(s)\right
){\frac {d}{ds}}y(s)}{{y}^{2}P-yP+Q{x}^{2}-Qx}}+\]\[+{ \frac
{\left ({\frac {\partial }{\partial y}}P(x,y)\right )y\left (y-1
\right )\left ({\frac {d}{ds}}y(s)\right
)^{2}}{{y}^{2}P(x,y)-yP(x,y)+ Q(x,y){x}^{2}-Q(x,y)x}}=0,
\end{equation}
 and
\[
\frac{ d^2 z}{ds^2}+A(s)\frac{d z}{ds}+B(s)\frac{d
t}{ds}+C(s)z(s)+E(s)t(s)=0,
\]
\[
\frac{ d^2 t}{ds^2}+F(s)\frac{d z}{ds}+H(s)\frac{d
t}{ds}+K(s)z(s)+L(s)t(s)=0
\]
with the coefficients depending from the parameter $s$ and the
functions $Q(x,y),P(x,y)$.

   Remark that last two equations are reduced at the independent equations
\[
\frac{ d^2 z}{ds^2}+M(s)\frac{d z}{ds}+N(s)z(s)=0
\]
and
\[
\frac{ d^2 t}{ds^2}+U(s)\frac{d t}{ds}+V(s)t(s)=0
\]
with the help of the first integral of geodesics
\[
z(s)\frac{ d x}{ds}+t(s)\frac{d y}{ds}-\alpha\frac{s}{2}-\beta=0
\]
of the metric (\ref{dryuma:eq21}).

\section{The second order ODE's cubic  on the first derivative in theory of the planar systems}

     The first two equations of geodesic of the metric (\ref{dryuma:eq23}) are
     equivalent to the one second order differential equation
\begin{equation} \label{dryuma:eq24}
{\frac {d^{2}}{d{x}^{2}}}y(x)+{\frac {\left (\left ({\frac
{\partial } {\partial y}}P(x,y)\right ){x}^{2}-\left ({\frac
{\partial }{\partial y}}P(x,y)\right )x\right )\left ({\frac
{d}{dx}}y(x)\right )^{3}}{{y}^
{2}P(x,y)-yP(x,y)+Q(x,y){x}^{2}-Q(x,y)x}}\!+\!\]\[\!+\!{\frac
{\left (\left ({ \frac {\partial }{\partial x}}P\!-\!{\frac
{\partial }{\partial y}}Q\right ){x}^{2}\!+\!\left ({\frac
{\partial }{\partial y}}Q\!-\!{ \frac {\partial }{\partial
x}}P\!-\!2\,P\right )x\!+\!\left ({\frac {\partial }{\partial
y}}P\right ){y}^{2}\!+\!\left (\!-\!2\,P\!-\!{ \frac {\partial
}{\partial y}}P\right )y\!+\!2\,P\right )\left ( {\frac
{d}{dx}}y(x)\right )^{2}}{{y}^{2}P\!-\!yP\!+\!Q{x}^{2}\!-\!Q
(x,y)x}}\!+\!\]\[+{\frac {\left (\!-\!\left ({\frac {\partial
}{\partial x}}Q \right ){x}^{2}\!+\!\left ({\frac {\partial
}{\partial x}}Q\!+\!2\,Q \right )x\!+\!\left ({\frac {\partial
}{\partial x}}P\!-\!{\frac {
\partial }{\partial y}}Q\right ){y}^{2}\!+\!\left (2\,Q\!-\!{\frac {
\partial }{\partial x}}P\!+\!{\frac {\partial }{\partial y}}Q
\right )y\!-\!2\,Q\right ){\frac {d}{dx}}y(x)}{{y}^{2}P-yP+
Q{x}^{2}-Qx}}\!+\!\]\[+{\frac {-\left ({\frac {\partial }{\partial
x} }Q(x,y)\right ){y}^{2}+\left ({\frac {\partial }{\partial
x}}Q(x,y) \right
)y}{{y}^{2}P(x,y)-yP(x,y)+Q(x,y){x}^{2}-Q(x,y)x}}=0.
\end{equation}

   The equation (\ref{dryuma:eq24}) has the equation
\[
\frac{d y}{dx}=\frac{Q(x,y)}{P(x,y)}
\]
as particular integral and the function
\begin{equation}\label{dryuma:eq25}
 y(x)=\frac{C(x-1)}{x-C}
\end{equation}
as the first integral.

    In result we get two-parametric solutions
of the second order ODE from the one-parametric solutions of a given first order
equation.

    This fact allow us to describe the properties of the first order system (\ref{dryuma:eq2})
with respect the values of parameter $C$ of the solution
(\ref{dryuma:eq25}).

    Let us consider some example.

    After substitution of the relation (\ref{dryuma:eq25}) into the
equation (\ref{dryuma:eq24}) one gets the expression
\begin{equation}\label{dryuma:eq26}
\alpha(x,y)C^5+\beta(x,y)C^4+\gamma(x,y)C^3+\delta(x,y)C^2+\epsilon(x,y)C+\mu(x,y)=0
\end{equation}, where
\[
\alpha(x,y)=\left ({\it a12}+{\it b12}\right ){y}^{2}+\]\[+\left
({\it b1}+{\it b2}+ \left (2\,{\it a22}+2\,{\it a11}+2\,{\it
b22}+2\,{\it b11}\right )x+{ \it a2}+{\it a1}\right )y+\left ({\it
a12}+{\it b12}\right ){x}^{2}+ \left ({\it a1}+{\it b2}+{\it
b1}+{\it a2}\right )x+\]\[+2\,{\it a0}+2\,{ \it b0}  ,
\]
\\[1mm]
\[
\beta(x,y)=2\,{\it b22}\,{y}^{3}+\left (\left (-5\,{\it a12}-{\it
b12}-2\,{\it b22}\right )x+2\,{\it b2}-{\it b12}\right
){y}^{2}+\]\[+\left (\left (\!-\!4\,{ \it a22}\!-\!3\,{\it
b2}\!-\!5\,{\it a1}\!-\!3\,{\it a2}\!-\!2\,{\it b11}\!-\!6\,{\it
b22}\!-\!{\it b1}\right )x\!-\!{\it a2}\right)+\]\[+\left(\!\left
(\!-\!4\,{\it b11}\!+\!{\it a12}\!-\!10\,{\it a11}\!-\!{\it
b12}\!-\!4\,{\it a22}\right ){x}^{2}\!+\!2\,{\it b0}\!-\!2\,{\it
b2}\!-\!{ \it b1}\right )y\!+\]\[+\!\left (\!-\!2\,{\it
a12}\!+\!2\,{\it a11}\right ){x}^{3}\!+\! \left (\!-\!2\,{\it
a12}\!-\!2\,{\it b1}\!-\!2\,{\it a1}\!-\!3\,{\it b12}\!-\!2\,{\it
a2} \right ){x}^{2}\!-\!4\,{\it b0}\!+\!\left (\!-\!6\,{\it
a0}\!-\!2\,{\it b1}\!-\!4\,{\it b0 }\!-\!3\,{\it b2}\!-\!{\it
a1}\!-\!2\,{\it a2}\right )x\!-\]\[-\!2\,{\it a0} ,
\]
\\[1mm]
\[
\gamma(x,y)=\left (-4\,x{\it b22}-2\,{\it b22}\right
){y}^{3}+\left (\left (-{\it b12}+2\,{\it b22}+10\,{\it a12}\right
){x}^{2}-2\,{\it b2}+\left ({ \it b12}+4\,{\it b22}-4\,{\it
b2}\right )x\right ){y}^{2}+\]\[+\left ({\it b2}\!+\!\left
(\!-\!3\,{\it a12}\!+\!2\,{\it a22}\!+\!2\,{\it b11}\!+\!{\it
b12}\!+\!20\,{\it a11}\right ){x}^{3}\!+\!\left (6\,{\it
b22}\!+\!2\,{\it a22}\!+\!{\it b1}\!+\!3\,{\it a2}\!-\!4\,{\it
b0}\!+\!6\,{\it b2}\right )x\right)y\!+\]\[+\left(\!\left
(\!-\!{\it b1}\!+\!10\,{\it a1}\!+\!3\, {\it a2}\!+\!2\,{\it
b2}\!+\!4\,{\it b11}\!+\!2\,{\it b12}\!-\!{\it a12}\!+\!8\,{\it
a22} \right ){x}^{2}\!-\!2\,{\it b0}\right )y\!+\!\]\[+\left ({\it
a12}-6\,{\it a11} \right ){x}^{4}+\left (4\,{\it a12}+{\it
a2}-2\,{\it a11}+{\it b1} \right ){x}^{3}+\left ({\it b1}+6\,{\it
a0}+3\,{\it b2}+{\it a2}+8\,{ \it b0}\right )x+\]\[+\left ({\it
a12}+6\,{\it a0}+4\,{\it b1}+2\,{\it a1}+ 2\,{\it b0}+3\,{\it
b12}+4\,{\it a2}\right ){x}^{2}+2\,{\it b0} ,
\]
\\[1mm]
\[
\delta(x,y)=\left (2\,{x}^{2}{\it b22}+4\,x{\it b22}\right
){y}^{3}+\left (\left ( {\it b12}-10\,{\it a12}\right
){x}^{3}+\left (2\,{\it b2}-4\,{\it b22} +{\it b12}\right
){x}^{2}+\left (4\,{\it b2}-2\,{\it b22}\right )x \right
){y}^{2}+\]\[+\left (\left ({\it b1}-4\,{\it a22}+3\,{\it
a12}-2\,{ \it b11}-{\it a2}-10\,{\it a1}-2\,{\it b12}\right
){x}^{3}+\left (-{ \it b12}-4\,{\it b2}-4\,{\it a22}-3\,{\it
a2}+{\it b1}+2\,{\it b0} \right ){x}^{2}\right)y+\]\[+\left(\left
(3\,{\it a12}-20\,{\it a11}\right ){x}^{4}+ \left (-3\,{\it
b2}+4\,{\it b0}-2\,{\it b22}\right )x\right )y+\]\[+6\,{x}^
{5}{\it a11}+\left (6\,{\it a11}+2\,{\it a1}-2\,{\it a12}\right
){x}^{ 4}+\left (-4\,{\it b0}-{\it b12}-2\,{\it b1}-2\,{\it
a2}-6\,{\it a0} \right ){x}^{2}+\]\[+\left (-2\,{\it a12}-2\,{\it
a0}-2\,{\it a2}-2\,{\it b1}\right ){x}^{3}+\left (-4\,{\it
b0}-{\it b2}\right )x ,
\]
\\[1mm]
\[
\epsilon(x,y)=-2\,{y}^{3}{x}^{2}{\it b22}+\left (-{x}^{3}{\it
b12}+5\,{x}^{4}{\it a12}+\left (-2\,{\it b2}+2\,{\it b22}\right
){x}^{2}\right ){y}^{2}+ \]\[+\left (\left (5\,{\it a1}-3\,{\it
a12}\right ){x}^{4}+\left (10\,{\it a11}-{\it a12}\right
){x}^{5}+\left (-{\it b1}+{\it b12}+{\it a2}+2\,{ \it a22}\right
){x}^{3}+\left (2\,{\it b2}-2\,{\it b0}\right ){x}^{2} \right
)y-2\,{x}^{6}{\it a11}+\]\[+\left (-{\it a1}-6\,{\it a11}\right
){x} ^{5}+\left ({\it b1}+2\,{\it a0}+{\it a2}\right
){x}^{3}+\left ({\it a12}-2\,{\it a1}\right
){x}^{4}+2\,{x}^{2}{\it b0} ,
\]
\\[1mm]
\[
\mu(x,y)=-{y}^{2}{x}^{5}{\it a12}+\left (-2\,{x}^{6}{\it
a11}+\left ({\it a12}- {\it a1}\right ){x}^{5}\right
)y+2\,{x}^{6}{\it a11}+{x}^{5}{\it a1}.
\]

      The substitution of the expression (\ref{dryuma:eq25}) into the relation (\ref{dryuma:eq26})
      give us the expression
\begin{equation}\label{dryuma:eq27}
 A(x)C^6+B(x)C^5+E(x)C^4+F(x)C^3+H(x)C^2+K(x)C+L(x)=0,
 \end{equation}
 where
\[
A(x)=\left (2\,{\it b12}-2\,{\it b22}-2\,{\it b11}+2\,{\it
a12}-2\,{\it a22 }-2\,{\it a11}\right ){x}^{2}+\]\[+\left (2\,{\it
b22}-2\,{\it b12}-2\,{ \it a12}+2\,{\it a11}+2\,{\it a22}+2\,{\it
b11}\right )x+\]\[+{\it a2}+{ \it a12}+{\it b1}+{\it b12}+{\it
b2}+2\,{\it a0}+{\it a1}+2\,{\it b0},
\]
\\[1mm]
\[
B(x)=\left (-8\,{\it a12}+12\,{\it a11}+4\,{\it b11}-4\,{\it
b22}+4\,{\it a22}\right ){x}^{3}+\]\[+\left (-2\,{\it b1}+7\,{\it
a12}+16\,{\it b22}+2\, {\it a1}+4\,{\it b2}-2\,{\it b11}-10\,{\it
a11}-5\,{\it b12}\right ){x }^{2}+\]\[+\left (-2\,{\it b1}-8\,{\it
b2}-8\,{\it a0}-4\,{\it a22}-4\,{ \it a2}+2\,{\it b12}-4\,{\it
a12}-14\,{\it b22}-8\,{\it b0}-6\,{\it a1 }-2\,{\it b11}\right
)x-\]\[-2\,{\it b0}-2\,{\it a0}-{\it b1}-{\it a2}-{ \it
b12}+2\,{\it b22}
\]
\\[1mm]
\[
E(x)=\left (-28\,{\it a11}-2\,{\it b11}-2\,{\it a22}-2\,{\it
b12}+12\,{\it a12}\right ){x}^{4}+\]\[+\left (-8\,{\it
a1}+18\,{\it a11}+6\,{\it b22}+4\, {\it b1}-2\,{\it b11}+4\,{\it
b12}-6\,{\it a22}-8\,{\it a12}-4\,{\it b2}\right
){x}^{3}+\]\[+\left (4\,{\it b1}+12\,{\it a0}+6\,{\it a2}+13\,{
\it a1}+{\it b2}+6\,{\it a12}+6\,{\it a22}+4\,{\it b11}+4\,{\it
b12}- 20\,{\it b22}+10\,{\it b0}\right ){x}^{2}+\]\[+\left
(4\,{\it a2}+2\,{\it b1}+10\,{\it b0}+10\,{\it b2}+2\,{\it
a22}+16\,{\it b22}+8\,{\it a0} \right )x-{\it b2}-2\,{\it b22},
\]
\\[1mm]
\[
F(x)=\left (32\,{\it a11}-8\,{\it a12}\right ){x}^{5}+\left
(-2\,{\it b1}+4 \,{\it a22}+12\,{\it a1}+3\,{\it b12}+2\,{\it
b11}+2\,{\it a12}-12\,{ \it a11}\right ){x}^{4}+\]\[+\left
(-4\,{\it a2}-6\,{\it b1}-12\,{\it a1}-4 \,{\it a12}-6\,{\it
b12}-8\,{\it a0}-2\,{\it b22}-2\,{\it b11}-4\,{ \it b0}+6\,{\it
b2}\right ){x}^{3}+\]\[+\left (-2\,{\it b1}-8\,{\it b2}-14 \,{\it
b0}-12\,{\it a0}-{\it b12}-6\,{\it a2}-4\,{\it a22}+6\,{\it b22
}\right ){x}^{2}+\left (-4\,{\it b22}-2\,{\it b0}-2\,{\it
b2}\right )x ,
\]
\\[1mm]
\[
H(x)=\left (2\,{\it a12}-18\,{\it a11}\right ){x}^{6}+\left
(2\,{\it a12}-2 \,{\it a11}-8\,{\it a1}\right ){x}^{5}+\]\[+\left
(-2\,{\it a22}+2\,{\it a0 }-{\it b12}+{\it a2}+3\,{\it b1}+3\,{\it
a1}+{\it a12}\right ){x}^{4}+ \]\[+\left (2\,{\it b1}+2\,{\it
b12}+6\,{\it b0}+4\,{\it a2}-2\,{\it b2}+8 \,{\it a0}+2\,{\it
a22}\right ){x}^{3}+\left (3\,{\it b2}+4\,{\it b0} \right
){x}^{2},
\]
\\[1mm]
\[
K(x)=4\,{x}^{7}{\it a11}+\left (6\,{\it a11}+2\,{\it a1}-{\it
a12}\right ){ x}^{6}+2\,{x}^{5}{\it a1}+\left (-2\,{\it a0}-{\it
a2}-{\it b1}\right ){x}^{4}-2\,{x}^{3}{\it b0},
\]
\\[1mm]
\[
L(x)=-2\,{x}^{7}{\it a11}-{x}^{6}{\it a1}.
\]

    The substitutions $x=0$ and $x=1$ into the (\ref{dryuma:eq27})
    lead to the conditions on the value $C$
\begin{equation} \label{dryuma:eq28}
\left ({\it a2}+{\it a12}+{\it b1}+{\it b12}+{\it b2}+2 \,{\it
a0}+{\it a1}+2\,{\it b0}\right ){C}^{2}+\]\[+\left (-2\,{\it
b0}-2\, {\it a0}-{\it b1}-{\it a2}-{\it b12}+2\,{\it b22}\right
)C-{\it b2}-2 \,{\it b22}=0,
\end{equation}
\\[1mm]
\begin{equation} \label{dryuma:eq29}
\left ({\it a2}+{\it a12}+{\it b1}+{\it b12}+{\it b2}+2\,{\it
a0}+{\it a1}+2\,{\it b0}\right ){C}^{2}+\]\[+\left (2 \,{\it
a11}-2\,{\it b0}-{\it b1}-2\,{\it a0}-{\it a12}-{\it a2}\right
)C-2\,{\it a11}-{\it a1}=0.
\end{equation}

   The substitution $(x=C, y=1-C)$ into the (\ref{dryuma:eq26}) lead
   to the conditions on the value $C$
\begin{equation} \label{dryuma:eq30}
\left ({\it b12}-2\,{\it b22}\right )C+2\,{\it b22}+{\it b2}=0.
\end{equation}

   After substitution $(y=1-x,x=1,C=1/C1)$ we get
\begin{equation} \label{dryuma:eq31}
\left ({\it a1}+2\,{\it a11}\right ){{\it C1}}^{2}+\left ({\it
a12}-2 \,{\it a11}+{\it b1}+{\it a2}+2\,{\it a0}+2\,{\it b0}\right
){\it C1}- {\it a1}-{\it a2}-{\it a12}-2\,{\it b0}-\]\[-{\it
b12}-{\it b1}-{\it b2}-2 \,{\it a0}=0,
\end{equation}
and the substitution $(y=1-x,x=0, C=1/C1)$ lead to the
condition
\begin{equation} \label{dryuma:eq32}
\left (-2\,{\it b22}-{\it b2}\right ){{\it C1}}^{2}+\left
(-2\,{\it b0 }-{\it b12}-{\it b1}-{\it a2}-2\,{\it a0}+2\,{\it
b22}\right ){\it C1} +{\it a1}+{\it a2}+{\it a12}+2\,{\it
b0}+\]\[+{\it b12}+{\it b1}+{\it b2}+2 \,{\it a0} =0.
\end{equation}

    Remark that for the comparison with results of the article (\cite{dryuma4:dryuma})
    in all above formulas the equation (\ref{dryuma:eq2}) was presented in the form
\[
\frac{dy}{dx}=\frac{a0+a1x+a2y+a11x^2+a12xy+a22y^2}{b0+b1x+b2y+b11x^2+b12xy+b22y^2}.
\]

\begin{rem}

   In the famous article (\cite{dryuma4:dryuma})
was developed the approach to the studying of the problem of the limit cycles of the
equations (\ref{dryuma:eq2}).

  Let us remind the basic facts of the Petrovsky-Landis theory.

   For every closed curve of the system (\ref{dryuma:eq2})  the solution
\[
 y(x)=\frac{C(x-1)}{(x-C)}
\]
 of equation
\begin{equation} \label{dryuma:eq33}
 \frac{d y}{dx}=-\frac{y(y-1)}{x(x-1)}
\end{equation}
is corresponded.

    At the same time the value $C$ satisfies the algebraic equations
\[
\sum a_n(\mu_i)C^n=0,
\]
where the coefficients $a_n(\mu_i)$ are dependent from the
parameters of equation (\ref{dryuma:eq2}).

    This equation arises  from the  condition

\[
\int_{c}{\frac{(x-C)^2[x(x-1)Q(x,y)+y(y-1)P(x,y)]}{x^3(x-1)^3}}dx=0.
\]

     The substitution of the function $y$ from the (\ref{dryuma:eq33}) into
     this expression and calculation of the residues with regard
     the points $0,1,C$ lead to the equations on the parameter
     $C$.

     According the (\cite{dryuma4:dryuma}) general quantity of the
     values $C$ defined by such equations is equal 14 and this
     number coincide with the quantity of closed solutions determined by the
     equations (\ref{dryuma:eq2}).

     As it was shown in (\cite{dryuma4:dryuma}) 11 curves from 14 can be transformed into the vicinity
     of the particular points of the equation (\ref{dryuma:eq30}).

     As result only three closed curves do not be transformed into  the vicinity
     of the particular points and the quantity of
     the limit cycles defined by the equation (\ref{dryuma:eq2}) was found  equal three.

       In spite of the fact that the  statement of the article (\cite{dryuma4:dryuma}) about
       the quantity of the limit cycles  defined by the system (\ref{dryuma:eq2}) was
       fallacious its approach is useful from point of our consideration.

       In fact the conditions on the value $C$ given by the (\ref{dryuma:eq28})- (\ref{dryuma:eq32})
       are the same with those which was used in the article (\cite{dryuma4:dryuma}) for
       estimation of the quantity of the limit cycles in quadratic
       polynomial system.

       It can be shown that all conditions on residues of functions
       considered in the article (\cite{dryuma4:dryuma}) are
       followed from the (\ref{dryuma:eq28})- (\ref{dryuma:eq30}).

\end{rem}

\section{The examples}

      We apply the second order ODE (\ref{dryuma:eq24}) for the studying of the properties of the
      first order ODE's
\[
 \frac{d y}{dx}=\frac{Q(x,y)}{P(x,y)}.
\]

   1. Let us consider the system
\begin{equation} \label{dryuma:eq34}
\dot y=8-3a-14ax-2axy-8y^2,
\]
\[
\dot x=2+4x-4ax^2+12xy.
\end{equation}

   Corresponding equation (\ref{dryuma:eq24}) for this system
   becomes
\begin{equation} \label{dryuma:eq35}
 \left
(\!-\!12y(x)^{3}x\!+\!\left
(\!-\!2\!+\!4\,a{x}^{2}\!+\!8\,{x}^{2}\right ) y(x)^{2}\!+\!\left
(2\,a{x}^{3}\!+\!2\!-\!6\,a{x}^{2}\!+\!4\,x\right )y
(x)\!-\!3\,ax\!-\!8\,{x}^{2}\!+\!8\,x\!-\!11\,a{x}^{2}\!+\!14\,a{x}^{3}
\right ){\frac {d^{2}}{d{x}^{2}}}y(x)+\]\[\left
(-12\,{x}^{3}+12\,{x}^{2}\right )\left ({\frac {d}{dx}}y(x) \right
)^{3}+\]\[+\left
(\!-\!4\,y(x){x}^{2}\!-\!4+4\,{x}^{2}+2\,a{x}^{2}+24\,y(x)x
-8\,a{x}^{2}y(x)-2\,a{x}^{3}+12\,\left (y(x)\right )^{2}x+4\,y(x)
\right )\left ({\frac {d}{dx}}y(x)\right )^{2}+\]\[\!+\!\left
(\!-\!12y(x)^{3}\!+\!\left (10\,ax\!+\!16\,x\!+\!8\right
)y(x)^{2}\!+\!\left (6\,a\!+\!20\,ax\!+\!2\,a{x}^{2}\!-\!12\right
)y(x)\!+\!16\!-\!8\,ax\!+\!14 \,a{x}^{2}\!-\!6\,a\!-\!16\,x \right
){\frac {d}{dx}}y(x)+\]\[+14\,ay(x )-12\,a\left (y(x)\right
)^{2}-2\,a\left (y(x)\right )^{3}=0
\end{equation}
 and the function
\begin{equation} \label{dryuma:eq36}
\frac{dy(x)}{dx}=\frac{8-3a-14ax-2axy-8y^2}{2+4x-4ax^2+12xy}.
\end{equation}
is solution of this equation.

  Moreover the function
\[
\frac{dy(x)}{dx}=-\frac{y(y-1)}{x(x-1)}
\]
 also is solution of the equation (\ref{dryuma:eq35}).

     Now we shall get another particular solutions of our equation (\ref{dryuma:eq35}).

     The substitution
\[
 \frac{d y}{dx}=\frac{M(x,y)}{N(x,y)}
\]
into the equation (\ref{dryuma:eq35}) give us the relation between
the function $M(x,y)$, $N(x,y)$.

     In spite of the fact that such relation is very complicated
     one can get with help of the MAPLE-6  some solutions of them.

     For example, the presentation of the functions $M(x,y)$,
     $N(x,y)$ in form
\[
M(x,y)=A(x)+B(x)y+C(x)y^2,\quad  N(x,y)=E(x)+F(x)y+H(x)y^2,
\]
where
\[
A(x)=2x-3ax^2-1,\quad B(x)=-1,\quad C(x)=-2x,
\]
\[
E(x)=x,\quad F(x)=2x^2,\quad H(x)=0,
\]
corresponds the system
\[
\frac{d y(x)}{dx}=\frac{2x-3ax^2-1-y-2xy^2}{x+2x^2y},
\]
compatible with the equation
\[
{\frac {d}{dx}}y(x)={\frac
{8-3\,a-14\,ax-2\,axy(x)-8\,y(x)^{2}}{2+4\,x-4\,a{x}^{2}+12\,y(x)x}}.
\]

   In  result we get the solution
\[
1/4+x-{x}^{2}+a{x}^{3}+xy(x)+{x}^{2}{y(x)}^{2}=0
\]
which present the limit cycle at the some value of the parameter
$a$.\\[2mm]

     2. Let us consider the system
\begin{equation} \label{dryuma:eq37}
\dot y=x+2y+4xy+(2+3a)y^2,
\]
\[
\dot x=5x+6x^2+4(1+a)xy+ay^2.
\end{equation}

   Corresponding equation (\ref{dryuma:eq24}) for this system
   becomes
\begin{equation} \label{dryuma:eq38}
{\frac {\left
(2\,{x}^{2}ay-4\,{x}^{2}a-2\,xay+4\,{x}^{3}a+4\,{x}^{3}-
4\,{x}^{2}\right )\left ({\frac {d}{dx}}y(x)\right
)^{3}}{\Delta(x,y)} }+\]\[+{\frac {\left
(-6\,xy-20\,y{x}^{2}-6\,x{y}^{2}a-3\,{x}^{2}-10\,{x}^{
2}ay-4\,{y}^{2}x-4\,{x}^{3}+7\,x+6\,xay\right )\left ({\frac
{d}{dx}}y (x)\right )^{2}}{\Delta(x,y)}}+\]\[+{\frac {\left
(6\,x{y}^{2}a-4\,a{y}^{2}
-7\,y+3\,{y}^{2}+20\,{y}^{2}x-x+{x}^{2}+4\,{y}^{3}-6\,xy+4\,y{x}^{2}+4
\,a{y}^{3}\right ){\frac {d}{dx}}y(x)}{\Delta(x,y)}}+\]\[+{\frac
{d^{2}}{d{ x}^{2}}}y(x)+{\frac
{y+3\,{y}^{2}-4\,{y}^{3}}{\Delta(x,y)}} =0,
\end{equation}
where
\[
\Delta(x,y)={y}^{4}a\!+\!\left (\!-\!a\!+\!4\,x\!+\!4\,xa\right
){y}^{3}\!+\!\left (\!-\!7\,xa\!-\!x\!+\!3\,{x}^{2
}a\!+\!8\,{x}^{2}\right ){y}^{2}\!+\!\left
(\!-\!8\,{x}^{2}\!-\!7\,x\!+\!4\,{x}^{3}\right
)y\!-\!{x}^{2}\!+\!{x}^{3}.
\]

   The first order equation
\begin{equation} \label{dryuma:eq39}
\frac{dy(x)}{dx}=\frac{x+2y+4xy+(2+3a)y^2}{5x+6x^2+4(1+a)xy+ay^2}.
\end{equation}
 is particular integral of the second order equation (\ref{dryuma:eq38}).

    The function $y(x)$ defined by the relation
\begin{equation} \label{dryuma:eq40}
x^2+x^3+x^2y(x)+2axy(x)^2+2axy(x)^3+a^2y(x)^4=0
\end{equation}
 satisfies  both equations (\ref{dryuma:eq38}) and
(\ref{dryuma:eq39}).

   Beyond this point the first order equation
\[
\frac{dy(x)}{dx}=-\frac{2x+3x^2+2xy+2ay^2+2ay^3}{x^2+4axy+6axy^2+4a^2y^3}
\]
resulting from the (\ref{dryuma:eq40}) by differentiation also
satisfies the equations (\ref{dryuma:eq38}).

     The presence of the two particular first order  integrals for a given second order differential equation
      can be used for finding one them in evident form.

      In concerned case we use the integral (\ref{dryuma:eq39}) as
      known and seek for the expression
\[
\frac{dy(x)}{dx}=\frac{M(x,y)}{N(x,y}
\]
as the second integral with some functions $M(x,y)$ and $N(x,y)$ (
for example as  polynomial on variable $y$ with the coefficients
depending from the variable $x$).

After substitution both particular integrals into a given second
order ODE (\ref{dryuma:eq38}) the functions $M(x,y)$ and $N(x,y)$
can be find in evident form from the corresponding condition of
compatibility.
\\[2mm]

    3. The next example is the system (\cite{dryuma6:dryuma})
\begin{equation} \label{dryuma:eq41}
\dot y=15(1+a)y+3a(1+a)x^2-2(9+5a)xy+16y^2,
\]
\[
\dot x=6(1+a)x+2y-6(2+a)x^2+12xy.
\end{equation}

   Corresponding equation (\ref{dryuma:eq24}) for this system
   becomes
\begin{equation} \label{dryuma:eq42}
{\frac {\left (-2\,x+12\,{x}^{3}-10\,{x}^{2}\right )\left ({\frac
{d}{ dx}}y(x)\right )^{3}}{\Delta(x,y)}}+\]\[+{\frac {\left
(16\,xy\!-\!20\,y{x}^{2}
\!+\!18\,{x}^{3}\!-\!39\,{x}^{2}\!-\!2\,{y}^{2}\!-\!12\,{y}^{2}x\!+\!21\,x\!-\!12\,yxa\!+\!10\,{x}
^{3}a\!-\!31\,a{x}^{2}\!+\!12\,y{x}^{2}a\!+\!21\,xa\!+\!2\,y\right
)\left ({\frac {d}{ dx}}y(x)\right
)^{2}}{\Delta(x,y)}}\!+\!\]\[\!+\!{\frac {\left
(6\,y{x}^{2}{a}^{2}\!-\!
18\,y{x}^{2}\!-\!21\,y\!+\!21\,{y}^{2}a\!-\!10\,{y}^{2}x\!+\!12\,{y}^{3}\!-\!21\,ya\!-\!22\,{y
}^{2}xa\!+\!9\,{y}^{2}\!-\!4\,y{x}^{2}a\!+\!54\,xy\!+\!42\,yxa\right
){\frac {d}{dx}}y (x)}{\Delta(x,y)}}\!+\!\]\[\!+\!{\frac
{d^{2}}{d{x}^{2}}}y(x)\!+\!{\frac {18\,{y}^{3}\!-\!
18\,{y}^{2}\!+\!6\,x{a}^{2}y\!-\!10\,{y}^{2}a\!-\!6\,x{a}^{2}{y}^{2}\!+\!10\,{y}^{3}a\!-\!
6\,{y}^{2}xa\!+\!6\,yxa}{\Delta(x,y)}}  =0,
\end{equation}
where
\[
\Delta(x,y)=\left (2+12\,x\right ){y}^{3}+\left
(-22\,x+4\,{x}^{2}-2-6\,a{x}^{2}+6 \,xa\right ){y}^{2}+\]\[+\left
(-10\,{x}^{3}a-18\,{x}^{3}+31\,a{x}^{2}-21\,
xa+45\,{x}^{2}-21\,x\right
)y-3\,{x}^{3}a-3\,{x}^{3}{a}^{2}+3\,{x}^{4} a+3\,{x}^{4}{a}^{2}.
\]

\section{The spatial first order system of equations}

    Generalization of the Petrovsky-Landis approach to the spatial
    first order systems of equations may be realized by the following
    way.

\subsection{}
     Instead of the spatial first order system of
equations
\begin{equation} \label{dryuma:eq43}
\frac{d x}{ds}=P(x,y,z),\quad \frac{d y}{ds}=Q(x,y,z),\quad
\frac{d z}{ds}=R(x,y,z)
 \end{equation}
 we consider the Pfaff equation
\begin{equation} \label{dryuma:eq44}
P(x,y,z)dx+Q(x,y,z)dy+R(x,y,z)dz=0 \end{equation} for the
orthogonal trajectories of the system (\ref{dryuma:eq43}).

 The relation (\ref{dryuma:eq44})  can be
considered as the linear first integral of geodesics of the the
six-dimensional space in local coordinates $(x,y,z,U,V,W)$ with
the Riemann metric of the form
\begin{equation} \label{dryuma:eq45}
^{6}ds^2=-2(\Gamma^1_{11}U+\Gamma^2_{11}V+\Gamma^3_{11}W)dx^2-2(\Gamma^1_{22}U+\Gamma^2_{22}V+\Gamma^3_{22}W)dy^2-\]\[-
2(\Gamma^1_{33}U+\Gamma^2_{33}V+\Gamma^3_{33}W)dz^2-4(\Gamma^1_{12}U+\Gamma^2_{12}V+\Gamma^3_{12}W)dxdy-\]
\[-4(\Gamma^1_{13}U+\Gamma^2_{13}V+\Gamma^3_{13}W)dxdz-4(\Gamma^1_{23}U+\Gamma^2_{23}V+\Gamma^3_{23}W)dydz+\]\[+2dxdu+2dydV+2dzdW,
\end{equation}
where $\Gamma^k_{ij}=\Gamma^k_{ij}(x,y,z)$.

    One part of geodesics of the metric (\ref{dryuma:eq45}) has the form of
deodesics of the three-dimensional space with  the affine
connection $\Gamma^k_{ij}=\Gamma^k_{ij}(x,y,z)$
\[
\frac{d^2x}{ds^2}+\Gamma^1_{11}\left(\frac{dx}{ds}\right)^2
+\Gamma^1_{22}\left(\frac{dy}{ds}\right)^2+\Gamma^1_{33}\left(\frac{dz}{ds}\right)^2+2\Gamma^1_{12}\frac{dx}{ds}\frac{dy}{ds}
+2\Gamma^1_{13}\frac{dx}{ds}\frac{dz}{ds}+2\Gamma^1_{23}\frac{dy}{ds}\frac{dz}{ds}=0,
\]
\[
\frac{d^2y}{ds^2}+\Gamma^2_{11}\left(\frac{dx}{ds}\right)^2
+\Gamma^2_{22}\left(\frac{dy}{ds}\right)^2+\Gamma^2_{33}\left(\frac{dz}{ds}\right)^2+2\Gamma^2_{12}\frac{dx}{ds}\frac{dy}{ds}
+2\Gamma^2_{13}\frac{dx}{ds}\frac{dz}{ds}+2\Gamma^2_{23}\frac{dy}{ds}\frac{dz}{ds}=0,
\]
\[
\frac{d^2z}{ds^2}+\Gamma^3_{11}\left(\frac{dx}{ds}\right)^2
+\Gamma^3_{22}\left(\frac{dy}{ds}\right)^2+\Gamma^3_{33}\left(\frac{dz}{ds}\right)^2+2\Gamma^3_{12}\frac{dx}{ds}\frac{dy}{ds}
+2\Gamma^3_{13}\frac{dx}{ds}\frac{dz}{ds}+2\Gamma^3_{23}\frac{dy}{ds}\frac{dz}{ds}=0
\]
and another part of geodesics is the linear system of the second
order differential equations for the coordinates $U,V,W$.

     Remind that if the relation
\[
a_i\frac{dx^i}{ds}=0
\]
is the linear integral of geodesics
\[
\frac{d^2x^i}{ds^2}+\Gamma^i_{jk}\frac{dx^j}{ds}\frac{dx^k}{ds}=0
\]
then the conditions
\[
a_{i;j}+a_{j;i}=0
\]
must be satisfy.

   In the space with the
   metric (\ref{dryuma:eq45}) these conditions allow us to determine only six coefficients
   $\Gamma^k_{ij}(x,y,z)$ from the eighteen one.

   For determination of others coefficients $\Gamma^k_{ij}(x,y,z)$
   we use the auxiliary equations same with the Petrovsky-Landis
   theory
\[
y(y-1)dx+x(x-1)dy=0,\quad z(z-1)dx+x(x-1)dz=0.
\]

Every of these relations can be considered as the linear first
integral of geodesics of the metric (\ref{dryuma:eq45}) and this
allow us to determined all coefficients of connection
$\Gamma^k_{ij}(x,y,z)$.

    Some of them are
\[
\Gamma^{1}_{11}(x,y,z)=\]\[={\frac {\left ({\frac {
\partial }{\partial x}}P(x,y,z)\right )x\left (x-1\right )}{P(x,y,z){x
}^{2}-P(x,y,z)x-Q(x,y,z){y}^{2}+Q(x,y,z)y-R(x,y,z){z}^{2}+R(x,y,z)z}},
\]
\[
\Gamma^{2}_{11}(x,y,z)=\]\[=-{\frac {y\left (y-1\right ) {\frac
{\partial }{\partial x}}P(x,y,z)}{P(x,y,z){x}^{2}-P(x,y,z)x-Q(x
,y,z){y}^{2}+Q(x,y,z)y-R(x,y,z){z}^{2}+R(x,y,z)z}},
\]
\[
\Gamma^{3}_{11}(x,y,z)=\]\[=-{\frac {\left ({\frac {
\partial }{\partial x}}P(x,y,z)\right )z\left (z-1\right )}{P(x,y,z){x
}^{2}-P(x,y,z)x-Q(x,y,z){y}^{2}+Q(x,y,z)y-R(x,y,z){z}^{2}+R(x,y,z)z}},
\]
\[
\Gamma^{2}_{13}(x,y,z)=\]\[=1/2\,{\frac {y\left (y-1 \right )\left
({\frac {\partial }{\partial x}}R(x,y,z) \left (x-{x}^{2}\right
)\!+\!2R(x,y,z)(x-1+z)\!-\!{ \frac {\partial }{\partial
z}}P(x,y,z)({x}^{2}-x)\right )}{ \left
(P(x,y,z){x}^{2}-P(x,y,z)x-Q(x,y,z){y}^{2}+Q(x,y,z)y-R(x,y,z){z}
^{2}+R(x,y,z)z\right )x\left (x-1\right )}}
\]
\[
\Gamma^{1}_{22}(x,y,z)={\frac {\left ({\frac {
\partial }{\partial y}}Q(x,y,z)\right )x\left (x-1\right )}{P(x,y,z){x
}^{2}-P(x,y,z)x-Q(x,y,z){y}^{2}+Q(x,y,z)y-R(x,y,z){z}^{2}+R(x,y,z)z}},
\]
\[
\Gamma^{2}_{22}(x,y,z)=-{\frac {y\left (y-1\right ) {\frac
{\partial }{\partial y}}Q(x,y,z)}{P(x,y,z){x}^{2}-P(x,y,z)x-Q(x
,y,z){y}^{2}+Q(x,y,z)y-R(x,y,z){z}^{2}+R(x,y,z)z}},
\]
\[
\Gamma^{3}_{22}(x,y,z)=-{\frac {z\left (z-1\right ) {\frac
{\partial }{\partial y}}Q(x,y,z)}{P(x,y,z){x}^{2}-P(x,y,z)x-Q(x
,y,z){y}^{2}+Q(x,y,z)y-R(x,y,z){z}^{2}+R(x,y,z)z}},
\]
\[
\Gamma^{1}_{23}(x,y,z)=1/2\,{\frac {x\left (x-1 \right )\left
({\frac {\partial }{\partial y}}R(x,y,z)+{\frac {
\partial }{\partial z}}Q(x,y,z)\right )}{P(x,y,z){x}^{2}-P(x,y,z)x-Q(x
,y,z){y}^{2}+Q(x,y,z)y-R(x,y,z){z}^{2}+R(x,y,z)z}},
\]
\[
\Gamma^{2}_{23}(x,y,z)=\]\[=-1/2\,{\frac {\left ({\frac {\partial
}{\partial y}}R(x,y,z)+{\frac {\partial }{\partial z}}Q(x,y,
z)\right )y\left (y-1\right
)}{P(x,y,z){x}^{2}-P(x,y,z)x-Q(x,y,z){y}^{
2}+Q(x,y,z)y-R(x,y,z){z}^{2}+R(x,y,z)z}},
\]
\[
\Gamma^{3}_{23}(x,y,z)=\]\[=-1/2\,{\frac {z\left (\left ({\frac
{\partial }{\partial y}}R(x,y,z)\right )z-{\frac {\partial }{
\partial y}}R(x,y,z)+\left ({\frac {\partial }{\partial z}}Q(x,y,z)
\right )z-{\frac {\partial }{\partial z}}Q(x,y,z)\right
)}{P(x,y,z){x}
^{2}-P(x,y,z)x-Q(x,y,z){y}^{2}+Q(x,y,z)y-R(x,y,z){z}^{2}+R(x,y,z)z}},
\]
\[
\Gamma^{1}_{33}(x,y,z)={\frac {x\left (x-1\right ){ \frac
{\partial }{\partial z}}R(x,y,z)}{P(x,y,z){x}^{2}-P(x,y,z)x-Q(x,
y,z){y}^{2}+Q(x,y,z)y-R(x,y,z){z}^{2}+R(x,y,z)z}},
\]
\[
\Gamma^{2}_{33}(x,y,z)=-{\frac {\left ({\frac {
\partial }{\partial z}}R(x,y,z)\right )y\left (y-1\right )}{P(x,y,z){x
}^{2}-P(x,y,z)x-Q(x,y,z){y}^{2}+Q(x,y,z)y-R(x,y,z){z}^{2}+R(x,y,z)z}},
\]
\[
\Gamma^{3}_{33}(x,y,z)=-{\frac {z\left (z-1\right ) {\frac
{\partial }{\partial z}}R(x,y,z)}{P(x,y,z){x}^{2}-P(x,y,z)x-Q(x
,y,z){y}^{2}+Q(x,y,z)y-R(x,y,z){z}^{2}+R(x,y,z)z}}.
\]

    So the six-dimensional space with the metric (\ref{dryuma:eq45}) is suitable geometric object for the
    studying of the properties of the spatial first order system of equations (\ref{dryuma:eq43}).

\subsection{}

    Another approach for the studying of the spatial
system of equations
    is connected with consideration of the system
    (\ref{dryuma:eq43})
\begin{equation}\label{dryuma:eq46}
 Q(x,y,z)dx-P(x,y,z)dy=0,\quad R(x,y,z)d x-P(x,y,z)dz=0
\end{equation}
 together with the relation
\begin{equation}\label{dryuma:eq47}
xdx+ydy+zdz=0, \quad x^2+y^2+z^2=C
\end{equation}
( or another ones) as the linear integrals of geodesics of the six-dimensional
space.

   From these conditions the eighteen  coefficients
of connections $\Gamma^k_{ij}$ can be uniquely determined.

    In fact, we get the following equations
\[
{\frac {d^{2}}{d{s}^{2}}}x(s)+{\frac {\left (\left ({\frac
{\partial } {\partial x}}R(x,y,z)\right )z+P(x,y,z)+y{\frac
{\partial }{\partial x }}Q(x,y,z)\right )\left ({\frac
{d}{ds}}x(s)\right )^{2}}{yQ(x,y,z)+R(
x,y,z)z+P(x,y,z)x}}+\]\[+{\frac {\left (-\left ({\frac {\partial
}{
\partial x}}P(x,y,z)\right )y+\left ({\frac {\partial }{\partial y}}Q(
x,y,z)\right )y+\left ({\frac {\partial }{\partial
y}}R(x,y,z)\right ) z\right )\left ({\frac {d}{ds}}x(s)\right
){\frac {d}{ds}}y(s)}{yQ(x,y ,z)+R(x,y,z)z+P(x,y,z)x}}-\]\[-{\frac
{\left (\left ({\frac {\partial }{
\partial y}}P(x,y,z)\right )y-P(x,y,z)\right )\left ({\frac {d}{ds}}y(
s)\right )^{2}}{yQ(x,y,z)+R(x,y,z)z+P(x,y,z)x}}-\]\[-{\frac {\left
(\left ( {\frac {\partial }{\partial x}}P(x,y,z)\right )z-\left
({\frac {
\partial }{\partial z}}R(x,y,z)\right )z-\left ({\frac {\partial }{
\partial z}}Q(x,y,z)\right )y\right )\left ({\frac {d}{ds}}x(s)\right
){\frac {d}{ds}}z(s)}{yQ(x,y,z)+R(x,y,z)z+P(x,y,z)x}}-\]\[-{\frac
{\left ( \left ({\frac {\partial }{\partial y}}P(x,y,z)\right
)z+y{\frac {
\partial }{\partial z}}P(x,y,z)\right )\left ({\frac {d}{ds}}y(s)
\right ){\frac
{d}{ds}}z(s)}{yQ(x,y,z)+R(x,y,z)z+P(x,y,z)x}}-{\frac { \left
(\left ({\frac {\partial }{\partial y}}P(x,y,z)\right )y-P(x,y,z
)\right )\left ({\frac {d}{ds}}z(s)\right
)^{2}}{yQ(x,y,z)+R(x,y,z)z+P (x,y,z)x}}=0,
\]
and analogous equations for the $${\frac
{d^{2}}{d{s}^{2}}}y(s)+...=0$$ and
 $$
 {\frac
{d^{2}}{d{s}^{2}}}z(s)+...=0
$$.

   In result the six-dimensional space with the geodesics, having the
   linear integral in form
   (\ref{dryuma:eq46})-(\ref{dryuma:eq47}) has been constructed.

\section{ Projectivization of the planar system}

    We apply the geometric approach at the studying of the Pfaff equations
    (\ref{dryuma:eq44}) connected with a planar systems of equations
\begin{equation}\label{dryuma:eq48}
 \frac{dx}{ds}=p(x,y),\quad \frac{dy}{ds}=q(x,y)
\end{equation}
after their projectivization.

    In fact after the continuation of the system (\ref{dryuma:eq46}) on the projective plane
    we get the Pfaff equation
\[
-zQ(x,y,z)dx+zP(x,y,z)dy+(xQ(x,y,z)-yP(x,y,z))dz=0
\]
where $P(x,y,z),Q(x,y,z)$ are the homogeneous functions
constructed with the help of given functions $p(x,y),q(x,y)$.

    As example we consider the  equation
\[
\frac{dx}{ds}=\lambda x-y-10x^2+(5+\delta)xy+y^2,\quad
\frac{dy}{ds}=x+x^2+(\epsilon-25)xy
\]
having at least a four limit cycles at the some restriction on the
 parameters $\lambda,\delta, \epsilon$.

    After its projectivization we get
    the Pfaff equation
\[
\left
(-{y}^{2}x\delta-y\lambda\,xz+z{y}^{2}-5\,{y}^{2}x+z{x}^{2}+{x}^
{3}+{x}^{2}y\epsilon-15\,{x}^{2}y-{y}^{3}\right ){\it
dz}+\]\[+\left ({z}^{
2}\lambda\,x-{z}^{2}y+5\,zxy+zxy\delta+z{y}^{2}-10\,z{x}^{2}\right
){ \it dy}+\left (-zxy\epsilon+25\,zxy-{z}^{2}x-z{x}^{2}\right
){\it dx} =0.
\]

   For determination of the connections coefficients of
   corresponding six-dimensional space we add at this equation the relations
\[
 y(y-1)dx+x(x-1)dy=0,\quad z(z-1)dx+x(x-1)dz=0
\]
and consider all of them as the linear integrals of geodesics.

    In result we get a six-dimensional space whose geodesics
    have at least a four limit cycles.

\subsection{The Lorenz system}

    The next example is the  system of equations
\[
\frac{dx}{ds}=P(x,y,z),\quad \frac{dy}{ds}=Q(x,y,z),\quad
\frac{dz}{ds}=R(x,y,z)
\]
with some functions $P,Q,R$.

     We consider   the relations
\[
Q(x,y,z)dx-P(x,y,z)dy=0,\quad R(x,y,z)dx-P(x,y,z)dz=0
\]
as the linear first integrals of geodesics of corresponding
six-dimensional space.

     For determination of the connection coefficients we use the invariant conditions
\[
\Gamma^1_{11}+\Gamma^2_{12}+\Gamma^3_{13}=0,\quad
\Gamma^1_{12}+\Gamma^2_{22}+\Gamma^3_{23}=0,\quad
\Gamma^1_{13}+\Gamma^2_{23}+\Gamma^3_{33}=0,
\]
and the freedom in the choice of coordinate systems.

    We use a following normalization of coordinate system
\[
\Gamma^1_{11}=0,\quad \Gamma^1_{12}=-\Gamma^3_{23},\quad
\Gamma^2_{22}=0,\quad \Gamma^2_{12}=-\Gamma^3_{13},\quad
\Gamma^3_{33}=0 \quad \Gamma^1_{13}=-\Gamma^2_{23}.
\]

     In result all connection coefficients of the space can be determined uniquely and
      we get the expressions for the connection coefficients
\[
\Gamma^{3}_{23}(x,y,z)=-1/4\,{\frac {-Q{ \frac {\partial
}{\partial y}}P+R{\frac {\partial }{
\partial z}}P+P{\frac {\partial }{\partial z}}R-2
\,P{\frac {\partial }{\partial x}}P+P{\frac {
\partial }{\partial y}}Q}{P Q}},
\]
\[
\Gamma^{2}_{33}(x,y,z)=-{\frac {Q(x,y,z){\frac {
\partial }{\partial z}}P(x,y,z)}{P(x,y,z)R(x,y,z)}},
\]
\[
\Gamma^{3}_{12}(x,y,z)=1/4\,{\frac {-R Q{\frac {\partial
}{\partial y}}P+\left (R\right )^{2}{ \frac {\partial }{\partial
z}}P+PR{\frac {
\partial }{\partial z}}R-2\,Q P{\frac {\partial }{
\partial y}}R-2\,P R{\frac {\partial }{\partial x}
}P+PR{\frac {\partial }{\partial y}}Q}{Q\left (P\right )^{2}}},
\]
\[
\Gamma^{2}_{13}(x,y,z)=-1/4\,{\frac {2\,Q P{\frac {\partial
}{\partial x}}P-Q P{\frac {
\partial }{\partial y}}Q-\left (Q\right )^{2}{\frac {
\partial }{\partial y}}P+QR{\frac {\partial }{
\partial z}}P-Q P{\frac {\partial }{\partial z}}R+
2\,\left ({\frac {\partial }{\partial z}}Q\right )R P}{\left
(P\right )^{2}R}},
\]
\[
\Gamma^{1}_{22}(x,y,z)=-{\frac {{\frac {\partial }{
\partial y}}P(x,y,z)}{Q(x,y,z)}},
\]
\[
\Gamma^{3}_{22}(x,y,z)=-{\frac {R(x,y,z){\frac {
\partial }{\partial y}}P(x,y,z)}{P(x,y,z)Q(x,y,z)}},
\]
\[
\Gamma^{1}_{23}(x,y,z)=-1/4\,{\frac {-2\,P{ \frac {\partial
}{\partial x}}P+P{\frac {\partial }{
\partial y}}Q+Q{\frac {\partial }{\partial y}}P+R
{\frac {\partial }{\partial z}}P+P{\frac {
\partial }{\partial z}}R}{R Q}},
\]
\[
\Gamma^{1}_{33}(x,y,z)=-{\frac {{\frac {\partial }{
\partial z}}P(x,y,z)}{R(x,y,z)}},
\quad \Gamma^{2}_{11}(x,y,z)=-{\frac {{\frac {\partial }{
\partial x}}Q(x,y,z)}{P(x,y,z)}},
\]
\[
\Gamma^{3}_{11}(x,y,z)=-{\frac {{\frac {\partial }{
\partial x}}R(x,y,z)}{P(x,y,z)}},
\]
\[ \Gamma^{2}_{23}(x,y,z)=1/4\,{\frac {2\,P{ \frac {\partial
}{\partial x}}P-P{\frac {\partial }{
\partial y}}Q-Q{\frac {\partial }{\partial y}}P+R
{\frac {\partial }{\partial z}}P-P{\frac {
\partial }{\partial z}}R}{P R}},
\]
\[
\Gamma^{3}_{13}(x,y,z)=-1/4\,{\frac {-P{ \frac {\partial
}{\partial y}}Q-Q{\frac {\partial }{
\partial y}}P+R{\frac {\partial }{\partial z}}P+P
{\frac {\partial }{\partial z}}R}{\left (P\right )^{2}}}.
\]

  The corresponding metric is
\begin{equation}\label{dryuma:eq49}
^{6}ds^2=-2\left
(\Gamma^{2}_{11}(x,y,z)V+\Gamma^{3}_{11}(x,y,z)W\right
)dx^{2}-\]\[-4\left
(-\Gamma^{3}_{23}(x,y,z)U-\Gamma^{3}_{13}(x,y,z)V+\Gamma^{3}_{12}(x,y,z)W\right
)dx dy-\]\[-4\left
(-\Gamma^{2}_{23}(x,y,z)U+\Gamma^{2}_{13}(x,y,z)V+\Gamma^{3}_{13}(x,y,
z)W\right )dx dz-\]\[-2\left (\Gamma^{1}_{22}
(x,y,z)U+\Gamma^{3}_{22}(x,y,z)W\right )dy^{ 2}-\]\[-4\left
(\Gamma^{1}_{23}(x,y,z)U+\Gamma^{2}_{23}(x,y,z)V+\Gamma^{3}_{23}(x,y,
z)W\right)dydz-\]\[-2\left(\Gamma^{1}_{33}(x,y,z)U+\Gamma^{2}_{33}(x,y,z)V\right
)dz^{ 2}+2dxdU+2dydV+2dzdW.
\end{equation}

    Remark that the conditions
\[
a_{i;j;k}+R^{m}_{k i j}a_m=0
\]
for the vectors
\[
a_i=[P(x,y,z),-Q(x,y,z),0,0,0,0]
\]
and
\[
a_i=[P(x,y,z),0,-R(x,y,z),0,0,0]
\]
are obeyed for the metric (\ref{dryuma:eq49}).

    In particular case of the Lorenz system of equations
\[
\frac{dx}{ds}=\sigma(y-x),\quad \frac{dy}{ds}=rx-y-xz,\quad
\frac{dz}{ds}=xy-bz,
\]
we have
\[
P(x,y,z)=\sigma(y-x),\quad Q(x,y,z)=rx-y-xz,\quad R(x,y,z)=xy-bz,
\]
 where
$\sigma,r, b$ are the parameters and geodesics of the metric
(\ref{dryuma:eq49}) are equivalent to the expressions
\[
{\frac {d^{2}}{d{x}^{2}}}y(x)-{\frac {\sigma\,\left ({\frac
{d}{dx}}y( x)\right )^{3}}{-rx+y+zx}}+1/2\,{\frac {\left
(x+xb-2\,x\sigma+rx-zx+2 \,\sigma\,y-yb-2\,y\right )\sigma\,\left
({\frac {d}{dx}}y(x)\right )^ {2}{\frac
{d}{dx}}z(x)}{r{x}^{2}y-yz{x}^{2}-{y}^{2}x+xb{z}^{2}-zxrb+yb
z}}-\]\[-{\frac {\left
(x+xb-2\,x\sigma+rx-zx+2\,\sigma\,y-yb-2\,y\right ) \left ({\frac
{d}{dx}}y(x)\right ){\frac {d}{dx}}z(x)}{bzx-ybz-y{x}^{2
}+{y}^{2}x}}-\]\[-1/2\,{\frac {\left
(-yb+xb+zx-2\,x\sigma+x-rx+2\,\sigma\, y\right )\left ({\frac
{d}{dx}}y(x)\right )^{2}}{rxy-{y}^{2}-yzx-r{x}^
{2}+xy+z{x}^{2}}}+\]\[+1/2\,{\frac {\left
(-yb+xb+zx-x-rx+2\,y\right ){ \frac {d}{dx}}y(x)}{\left
({y}^{2}-2\,xy+{x}^{2}\right )\sigma}}-{ \frac {z-r}{\sigma\,\left
(-y+x\right )}} =0,
\]
\\[1mm]
\[
{\frac {d^{2}}{d{x}^{2}}}z(x)+1/2\,{\frac {\left
(x+xb-2\,x\sigma+rx-z x+2\,\sigma\,y-yb-2\,y\right )\sigma\,\left
({\frac {d}{dx}}y(x) \right )\left ({\frac {d}{dx}}z(x)\right
)^{2}}{r{x}^{2}y-yz{x}^{2}-{y
}^{2}x+xb{z}^{2}-zxrb+ybz}}-\]\[-{\frac {\sigma\,\left ({\frac
{d}{dx}}y(x) \right )^{2}{\frac {d}{dx}}z(x)}{-rx+y+zx}}-{\frac
{\left (-yb+xb+zx-2 \,x\sigma+x-rx+2\,\sigma\,y\right )\left
({\frac {d}{dx}}y(x)\right ){ \frac
{d}{dx}}z(x)}{rxy-{y}^{2}-yzx-r{x}^{2}+xy+z{x}^{2}}}+\]\[+{\frac {
\left (-xy+bz\right )\left ({\frac {d}{dx}}y(x)\right
)^{2}}{rxy-{y}^{ 2}-yzx-r{x}^{2}+xy+z{x}^{2}}}-1/2\,{\frac {\left
(x+xb-2\,x\sigma+rx-z x+2\,\sigma\,y-yb-2\,y\right )\left ({\frac
{d}{dx}}z(x)\right )^{2}}{
bzx-ybz-y{x}^{2}+{y}^{2}x}}+\]\[\!+\!1/2\,{\frac {\left
(xb{z}^{2}\!+\!\left ( \left
(\!-\!3\,{x}^{2}\!-\!{b}^{2}\!+\!2\,\sigma\,b\right
)y\!+\!2\,{x}^{3}\!+\!\left ({b}
^{2}\!-\!2\,\sigma\,b\!+\!b-\!rb\right )x\right )z\!+\!\left
(\!-\!2\,\sigma\!+\!b\!-\!2\right ) x{y}^{2}\!+\!\left
(3\,r\!-\!b\!+\!2\,\sigma\!+\!1\right
){x}^{2}y\!-\!2\,{x}^{3}r\right ) {\frac {d}{dx}}y(x)}{\left
(-y+x\right )^{2}\left (-rx+y+zx\right )
\sigma}}\!-\!\]\[\!-\!1/2\,{\frac {\left
(-yb+xb+zx-x-rx+2\,y\right ){\frac {d}{dx} }z(x)}{\left
({y}^{2}-2\,xy+{x}^{2}\right )\sigma}}+{\frac {y}{\left (
-y+x\right )\sigma}} =0,
\]
and the second order linear system differential equations for the
coordinates $U(s),V(s),W(s)$.

   These equations have  the solutions in form  of the first order differential
equations
\[
\frac{dy}{dx}=\frac{rx-y-xz}{\sigma(y-x)},\quad
\frac{dz}{dx}=\frac{xy-bz}{\sigma(y-x)}.
\]

     Geometrical properties of the metric (\ref{dryuma:eq49}) are dependent
     from the parameters $\sigma,b,r$.

     For example the component $R_{zz}$ the Ricci tensor $R_{ij}$
     look as
\[
R_{zz}= -1/2\,{\frac {x\left (-zx+\left (2\,\sigma-b-2\right
)y+\left (1+b-2\, \sigma+r\right )x\right )}{\left (-xy+bz\right
)\left (-y+x\right ) \left (-rx+y+zx\right )}}
\]
and a more complicated expressions for another components.

    Finally  we present the expression for the density of the  Chern-Simons invariant of affine
   connection for the space defined by the Lorenz system of
   equation.

   In case of three dimensional  space it is defined as  (\cite{dryuma6:dryuma})
\begin{equation} \label{dryuma:eq50}
CS(\Gamma)=\epsilon^{i j k}(\Gamma^p_{i q}\Gamma^q_{k
p;j}+\frac{2}{3}\Gamma^p_{i q}\Gamma^q_{j r}\Gamma^r_{k p})
\end{equation}
where $\epsilon^{i j k}$ is a Levi-Civita symbols.

   For the six-dimensional metric we get the expression
\begin{equation}\label{dryuma:eq51}
4\,\sigma\,\left (-y+x\right )^{4}\left (-rx+y+zx\right )^{2}\left
(-x y+bz\right )^{2} CS(\Gamma)=\]\[=-9\,{z}^{5}{x}^{4}b+\left
(\left (10\,{x}^{5}+\left (11\,{b}^{2}+48\, \sigma\,b-40\,b\right
){x}^{3}\right )y+\left (-11\,{b}^{2}+4\,b+36\,r
b-25\,\sigma\,b\right ){x}^{4}-23\,\sigma\,{y}^{2}{x}^{2}b-{x}^{6}
\right ){z}^{4}+\]\[+M3(x,y){z}^{3}+M4(x,y){z}^{2}+M5(x,y
)z+M6(x,y), \end{equation}
 where
\[
M3(x,y)=\left (4\,r-2\,\sigma\right ){x}^{6}+\left
(18\,b-40\,r-8+32\,\sigma \right )y{x}^{5}+\]\[+\left (\left
(-18\,b+44-56\,\sigma\right ){y}^{2}-14
\,\sigma\,{b}^{2}-54\,b{r}^{2}+12\,{b}^{2}+7\,b-18\,{\sigma}^{2}b-4\,
\sigma\,b+75\,\sigma\,rb+5\,{b}^{3}-12\,rb+26\,r{b}^{2}\right
){x}^{4} +\]\[+\left (26\,\sigma\,{y}^{3}+\left
(-144\,\sigma\,rb-2\,b+80\,{\sigma}^
{2}b+120\,rb-66\,\sigma\,b-10\,{b}^{3}-50\,{b}^{2}-26\,r{b}^{2}+14\,
\sigma\,{b}^{2}\right )y\right ){x}^{3}+\]\[+\left
(-106\,{\sigma}^{2}b+14
\,\sigma\,{b}^{2}+5\,{b}^{3}+38\,{b}^{2}+69\,\sigma\,rb-59\,b+138\,
\sigma\,b\right ){y}^{2}{x}^{2}+\left
(44\,{\sigma}^{2}b-68\,\sigma\,b -14\,\sigma\,{b}^{2}\right
){y}^{3}x.
\]

The functions $M4(x,y)-M6(x,y)$ are dependent from the parameters
and are too cumbersome.

   Remark that in the case $x=y$ for the right part of the (\ref{dryuma:eq51}) we
   get the expression
\[
-9\,{y}^{4}\left (z+1-r\right )^{4}\left (bz-{y}^{2}\right ).
\]
which is equal to zero on the stationary  points of the Lorenz model
\[
z=r-1,\quad y=\sqrt{b(r-1)}.
\]

    More detail information on the properties of the
    six-dimensional Riemann space connected with the Lorenz system of
    equations can be obtained with the help of the studying of the linear system of
    equations for the coordinates $U,V,W$.

    The studying of the equation
\[
g^{ij}\nabla_i\nabla_j A_k-R^{l}_{k}A_l+\lambda A_k=0
\]
for the eigenvalues $\lambda$  of the de Rham operator, defined on the 1-forms
\[
A(x,y,z)=A_i(x,y,z)dx^i
\]
of manifold can be also useful for that.

\end{document}